\documentclass{amsart}

\usepackage{amssymb}

\usepackage[all]{xy}

\newtheorem{thm}{Theorem}

\newtheorem{lem}[thm]{Lemma}
\newtheorem{cor}[thm]{Corollary}

\theoremstyle{definition}
\newtheorem{defn}[thm]{Definition}

\newtheorem{say}[thm]{}
\newtheorem{exmp}[thm]{Example}

\newtheorem*{ack}{Acknowledgments}      

\newtheorem{defn-thm}[thm]{Definition--Theorem}  
\newtheorem{defn-lem}[thm]{Definition--Lemma}  

\theoremstyle{remark}

\renewcommand{\c}[0]{{\mathbb C}}  

\renewcommand{\o}[0]{{\mathcal O}} 

\renewcommand{\a}[0]{{\mathbb A}}

\newcommand{\q}[0]{{\mathbb Q}}
\newcommand{\map}[0]{\dasharrow}
\newcommand{\qtq}[1]{\quad\mbox{#1}\quad}
\newcommand{\spec}[0]{\operatorname{Spec}}

\newcommand{\supp}[0]{\operatorname{Supp}}    
\newcommand{\red}[0]{\operatorname{red}}

\newcommand{\simq}[0]{\sim_{\q}}

\newcommand{\tr}[0]{\operatorname{tr}}

\newcommand{\mld}[0]{\operatorname{mld}}

\begin{document}

\bibliographystyle{amsalpha}


\title[Seminormal log centers]{Seminormal log centers and \\ 
deformations of pairs}
\author{J\'anos Koll\'ar}

\maketitle

The philosophy of Shokurov \cite{sho-3ff} stresses the
importance of understanding the log canonical centers 
 of an lc pair $(X,\Delta)$ (see Definition \ref{log.cent.defn}). 
After the initial work of
\cite{kawa98}, a systematic study was started by
 \cite{ambro}. For  extensions, surveys and comprehensive treatments
see  \cite{amb-surv} and   \cite{fujinobook}.
The following are  two of their principal results.
\begin{itemize}
\item Any union of log canonical centers is seminormal
(see Definition \ref{rel.semin.defn}).
\item Any intersection of  log canonical centers is
also a union of log canonical centers.
\end{itemize}
The aim of this note is to extend these results to
certain subvarieties of  an
lc pair $(X,\Delta)$ that are close to being a  log canonical center.
To state our results, we need a definition.
(See \cite{km-book} for basic concepts and results  related to MMP.
As in the above papers, we also work over a field of characteristic 0.)

\begin{defn}\label{log.cent.defn} Let $(X,\Delta)$ be lc
and $Z\subset X$ an irreducible subvariety. 
Following Shokurov and Ambro, the {\it minimal log discrepancy}
of $Z$ is the infimum of the numbers
$1+a(E,X,\Delta)$ as $E$ runs through all divisors over $X$
whose center is $Z$ \cite{amb-min}.
It is denoted by $\mld(Z,X,\Delta)$.

An irreducible subvariety 
$Z\subset X$ is called a {\it log center} of  $(X,\Delta)$
if  $\mld(Z,X,\Delta)<1$.
If $Z\subset X$ is a divisor, then $Z$ is a log center iff
it is an irreducible component of $\Delta$ and then
its coefficient is $1-\mld(Z,X,\Delta)$.

A {\it log canonical center}  is  a 
log center whose minimal log discrepancy equals  $0$.
\end{defn}

Our first aim is to prove the following.
(See Definition \ref{rel.semin.defn} for seminormality.)

\begin{thm} \label{1/6.logcent.semin.thm}
Let $(X,\Delta)$ be an lc pair and $Z_i\subset X$   log centers 
 for  $i=1,\dots, m$.
\begin{enumerate}
\item If  $\mld(Z_i,X,\Delta)<\frac16$ for every $i$ 
then $Z_1\cup\cdots\cup Z_m$ is
 seminormal.
\item If $ \sum_{i=1}^m\mld(Z_i,X,\Delta)<1$ then 
every irreducible component of $Z_1\cap\cdots\cap Z_m$
is a log center with minimal log discrepancy
$\leq \sum_{i=1}^m\mld(Z_i,X,\Delta)$.
\end{enumerate}
\end{thm}

A result of this type is not entirely surprising.
By Shokurov's conjecture on the boundedness of complements
 \cite{sho-3ff}, if $(X,\sum a_iD_i)$ is lc and the $a_i$ are close enough to 1,
then there is another lc pair  $(X,\Delta'+\sum D_i)$ where the $D_i$ all 
appear with coefficient 1. 
Thus the $D_i$ are  log canonical centers of $(X,\Delta'+\sum D_i)$
hence their union is seminormal and Du~Bois  \cite{k-db}.
In particular,
there should be  a function $\epsilon(n)>0$ such that
the union of the $D_i$ with $a_i> 1-\epsilon(\dim X)$  is
 seminormal and Du~Bois.
 The function  $\epsilon(n)$ is not known,
but it must converge to 0 at least doubly exponentially.
(See \cite[Sec.8]{k-pairs} for the  conjectured optimal 
value of $\epsilon(n)$ and
for examples.)

Thus it is somewhat unexpected that, at least for seminormality,
the bound in  Theorem \ref{1/6.logcent.semin.thm}.1
is independent of the dimension.

Note that  we do not assert that these $Z_i$ 
are log canonical centers of
 some other lc pair $(X,\Delta')$; this is actually not true.
In particular, unlike  log canonical centers, the $Z_i$ 
are not Du~Bois in general; see Example \ref{1/6.logcent.exmps}.5. 

As  Examples  \ref{1/6.logcent.exmps}.1--3 show, the value
$\frac16$ is optimal. There is, however, one important special case
when it can be improved to $\frac12$. The precise statement
 is given in Theorem \ref{1/2.logcent.semin.thm}; here we mention 
a consequence which was the main reason of this project.
The result implies that if we consider the moduli of
lc pairs  $(X,\Delta)$ where all the coefficients in $\Delta$
are $>\frac12$, then we do not have to worry about embedded
points on $\Delta$. (Examples of Hassett show that  embedded
points do appear when  the coefficients in $\Delta$
are $\leq \frac12$. See \cite[Sec.6]{k-modsurv} for an overview
and the forthcoming \cite{k-modbook} for details.)

\begin{cor} \label{1/2.coeff.seminorm.thm}
Let $(X,\Delta=\sum_{i\in I} b_i B_i)$ be lc. 
Let $f:X\to C$ be a morphism to a smooth curve
such that $(X,X_c+\Delta)$ is lc for every fiber $X_c:=f^{-1}(c)$.
 Let $J\subset I$ be any subset such that
$b_j>\frac12$ for every $j\in J$ and set $B_J:=\cup_{j\in J} B_j$. 

Then $B_J\to C$ is flat with reduced fibers. 
\end{cor}

The extension of these results to the semi log canonical
case requires additional considerations; these will be
treated in \cite[Chap 3]{k-modbook}.

\begin{say}\label{1/2.mmp.assumption}
The proof of  Theorem \ref{1/6.logcent.semin.thm} uses the following
 recently established result of Birkar \cite{bir} and Hacon and Xu \cite{hacon-xu}.
For $\dim X\leq 4$, it also follows from earlier results of
Shokurov  \cite{sho-VII}.
\medskip

{\it Theorem} \ref{1/2.mmp.assumption}.1. 
Let $g:X\to S$ be a projective, birational morphism and
 $\Delta',\Delta''$  effective $\q$-divisors on $X$ such that
 $\bigl(X, \Delta'+\Delta''\bigr)$ is dlt, $\q$-factorial and
 $K_X+\Delta'+\Delta''\sim_{\q,g}0$.
Then the $\bigl(X, \Delta''\bigr)$-MMP  with scaling over $S$ terminates with
 a $\q$-factorial minimal model.
\medskip

One of the difficulties in \cite{ambro, fujinobook} comes
from making the proof independent of MMP assumptions.
The proof in \cite{hacon-xu} uses several delicate
properties of log canonical centers, including some 
of the theorems of \cite{ambro, fujinobook}.
Thus, although the statement of Theorem \ref{1/6.logcent.semin.thm}
sharpens several of the  theorems of \cite{ambro, fujinobook}
on lc centers, it does not give a new proof.

\end{say}

\begin{exmp}\label{1/6.logcent.exmps}
The following examples show that the numerical conditions of
Theorem \ref{1/6.logcent.semin.thm} are sharp.

 1.  $\bigl(\a^2, \frac56(x^2=y^3)\bigr)$ is lc,
the curve $(x^2=y^3)$ is a log center with $\mld=\frac16$ but it is not
seminormal.

2. Consider  
$\bigl(\a^3, \frac{11}{12}(z-x^2-y^3)+\frac{11}{12}(z+x^2+y^3)\bigr)$.
One can check  that this is lc.
The irreducible components of the boundary are smooth, but
their intersection is a cuspidal curve, hence not seminormal.
It is again a log center with $\mld=\frac16$.

3. 
The image of $\c^2_{uv}$ by the map $x=u, y=v^3, z=v^2, t=uv$ 
is a divisor $D_1\subset X:=(xy-zt)\subset \c^4$
and $\c^2\to D_1$ is an isomorphism outside the origin.
Note that  the zero set of $(y^2-z^3)$ is
$D_1+2(y=z=0)$. Let $D_2, D_3$ be 2 general members
of the family of planes in the linear system $|(y=z=0)|$.
We claim that  $\bigl(X, \frac56 D_1+\frac56 D_2+\frac56 D_3\bigr)$ is lc.
Here $D_1$ is a log center  with $\mld=\frac16$ but  seminormality fails in 
codimension 3 on $X$.

In order to check the claim, blow up the ideal $(x,z)$.
On $\c^2_{uv}$ this corresponds to blowing up the ideal $(u,v^2)$.

On one of the charts we have coordinates
$x_1:=x/z, y,z$ and the birational transform $D'_1$ of
$D_1$ is given by $(y^2=z^3)$. 
On the other  chart we have coordinates
$x, z_1:=z/x,t$ and 
$D'_1$ is given by $(z_1x^3=t^2)$.  
Thus we see that $\bigl(B_{(x,z)}X, \frac56 D'_1\bigr)$ is lc.
The linear system  $|(y=z=0)|$ becomes base point free on the blow-up,
hence 
 $\bigl(B_{(x,z)}X, \frac56 D'_1+\frac56 D'_2+\frac56 D'_3\bigr)$ is lc
and so is $\bigl(X, \frac56 D_1+\frac56 D_2+\frac56 D_3\bigr)$.

4. Assume that $(X,\sum_{i\in I} a_iD_i)$ has simple normal crossing
and $a_i\leq 1$ for every $i$.
Let $J\subset I$ be a subset such that $a_j> 0$ for every $j\in J$
and $\sum_{j\in J}a_j> |J|-1$. Then every irreducible component of
$\cap_{j\in J}D_j$ is a  log center of  $(X,\sum_{i\in I} a_iD_i)$
with $\mld=\sum_{j\in J}(1-a_j)=|J|-\sum_{j\in J}a_j$.
In particular, $D_i$ is 
a    log center of  $(X,\sum_{i\in I} a_iD_i)$ with $\mld=1-a_i$.
Thus Theorem \ref{1/6.logcent.semin.thm}.2 is sharp.
  By \cite[Sec.2.3]{km-book}, every  log center of  $(X,\sum_{i\in I} a_iD_i)$
arises this way. 

5. Let $X$ be a smooth variety and $D\subset X$ a reduced divisor.
Then $D$ is Du~Bois iff $(X,D)$ is lc.
(See \cite{kss10, k-db} for much stronger results.)
In particular,  $D:=(x^2+y^3+z^7=0)\subset \a^3$
is a log center of 
the lc pair $\bigl(\a^3, \frac{41}{42}D\bigr)$ with $\mld= \frac1{42}$ but
$D$ is not Du~Bois and it can not be an lc center of any lc pair
$(X,\Delta)$. On the other hand, $D$ is normal hence seminormal.
\end{exmp}

\begin{say}[Log centers and birational maps]\label{log.cent.bir.say}
Let $g:(Y,\Delta_Y)\to (X,\Delta_X)$ be a proper birational morphism
between lc pairs (with $\Delta_X,\Delta_Y $ not necessarily effective) such that
$K_Y+\Delta_Y\simq g^*\bigl(K_X+\Delta_X\bigr)$
and $g_*\Delta_Y=\Delta_X$.

If $Z\subset Y$ is a log center of
$(Y,\Delta_Y)$ then $g(Z)$ is also a log center of
$(X,\Delta_X)$ with the same $\mld$. Moreover, every log center of
$(X,\Delta_X)$ is the image of a log center of
$(Y,\Delta_Y)$.

Thus, for any $(X,\Delta_X)$, we can use a log resolution
$g:(Y,\Delta_Y)\to (X,\Delta_X)$ to reduce the computation of
log centers to the simple normal crossing
 case considered in  Example \ref{1/6.logcent.exmps}.4.

This implies that an  lc pair $(X,\Delta)$ 
has  only finitely many log centers and 
the union of all log centers of codimension $\geq 2$
is the smallest closed
subscheme $W\subset X$ such that $\bigl(X\setminus W,\Delta|_{X\setminus W}\bigr)$
 is canonical.
\end{say}

\begin{say}[Proof of the divisorial case of Theorem \ref{1/6.logcent.semin.thm}]
\label{5/6.first.emp}
We show Theorem \ref{1/6.logcent.semin.thm}
in the  special case 
when $(X,\Delta')$ is dlt for some $\Delta'$ 
  and the $Z_i=:D_i$ are $\q$-Cartier divisors.

Since  $(X,\Delta')$ is dlt, $X$  is Cohen--Macaulay and so is $\sum D_i$
\cite[5.25]{km-book}. 
In particular, $\sum D_i$ satisfies Serre's condition  $S_2$.
An $S_2$-scheme is  seminormal iff it 
 is  seminormal at its codimension 1 points. 
By localization at  codimension 1 points we are reduced to the case
when $\dim X=2$. 

Then $X$ has a quotient singularity at every point of
$\sum D_i$, and Reid's covering method \cite[20.3]{k-etal} reduces the claim
to the smooth case.  It is now an elementary exercise
to see that if $(\a^2, \sum a_iD_i)$ is lc and $a_i>\frac 56$
then $\sum D_i$ has only ordinary nodes, hence
it is seminormal.

Next we  prove Theorem \ref{1/6.logcent.semin.thm}.2
 assuming  that $m=2$ and $Z_i=:D_i$
are  $\q$-Cartier divisors. Then every irreducible component
of $D_1\cap D_2$ has codimension 2, thus it is again 
enough to check the smooth  surface case. 
The exceptional divisor of the blow up of $x\in D_1\cap D_2$
 shows that $x$ is a log center with
$\mld\leq (1-a_1)+(1-a_2)$.
\end{say}

Any argument along this line breaks down completely
if we only assume that $(X, \sum a_iD_i)$ is lc.
In general the $D_i$ are not $S_2$, not even if $a_i=1$.
Thus seminormality at  codimension 1 points does not imply
seminormality.

Instead, we choose a suitable dlt model $(Y, \Delta_Y)$ of  $(X, \Delta)$,
use (\ref{5/6.first.emp}) on it and then
descend seminormality from $Y$  to $X$.
The next two lemmas will be used to construct $(Y, \Delta_Y)$.

\begin{lem}\label{dlt.model.withD.lem}
 Let $(X,\Delta)$ be lc. Then there is a 
projective, birational morphism 
$g:(Y,\Delta_Y)\to (X,\Delta)$ such that
\begin{enumerate}
\item $\bigl(Y,\Delta_Y\bigr)$ is dlt, $\q$-factorial
 (and $\Delta_Y$ is effective),
\item $K_Y+\Delta_Y\simq g^*\bigl(K_X+\Delta\bigr)$,
\item for every log center $Z\subset X$ of $(X,\Delta)$ there is a
divisor $D_Z\subset Y$ such that $g(D_Z)=Z$ and
 $D_Z$ appears in $\Delta_Y$ with coefficient $1-\mld(Z,X,\Delta) $.
\end{enumerate}
\end{lem}

Proof. This is well known.
Under suitable MMP assumptions, a proof is given in
\cite[17.10]{k-etal}. One can remove the  MMP assumptions
as follows. 

A method of Hacon (cf.\ \cite[3.1]{k-db})
constructs a model satisfying  (1--2). 
Since there are only finitely many log centers, 
it is enough to add the divisors $D_Z$ one at  a time.
This is explained in \cite[37]{k-acc}. 
A simplified proof is in \cite[Sec.4]{fuj-ssmmp}.
\qed

\begin{lem} \label{dlt.model.nefD.lem}
Let $g:Y\to X$ be a
projective, birational morphism  and
$\Delta_1,\Delta_2$  effective $\q$-divisors on $Y$.
Assume that 
 $\bigl(Y,\Delta_1+\Delta_2\bigr)$ is dlt, $\q$-factorial and
 $K_Y+\Delta_1+\Delta_2\sim_{\q,g} 0$. 
By (\ref{1/2.mmp.assumption}.1) a suitable  $(Y, \Delta_2)$-MMP   over $X$
terminates with a  $\q$-factorial minimal model
$g^m:\bigl(Y^m, \Delta_2^m\bigr)\to X$. Then
\begin{enumerate}
\item $-\Delta_1^m$ is $g^m$-nef,
\item  $g\bigl( \Delta_1\bigr)=g^m\bigl( \Delta_1^m\bigr)$ and
\item $\supp (g^m)^{-1}\bigl(g^m\bigl( \Delta_1^m\bigr)\bigr)=
\supp \Delta_1^m$.
\end{enumerate}
\end{lem}

Proof. Since $K_Y+\Delta_1+\Delta_2\sim_{\q,g} 0$,
we see that $K_{Y^m}+\Delta_1^m+\Delta_2^m\sim_{\q,g^m} 0$.
Thus $-\Delta_1^m \sim_{\q,g^m}K_{Y^m}+\Delta_2^m$ is $g^m$-nef.
Since $g^m$ has connected fibers and $\Delta_1^m$ is effective,
every fiber of $g^m$ is either contained in
$\supp \Delta_1^m$ or is disjoint from it.
This proves (3). 

In order to see (2), we prove by induction
that, at every intermediate step $g^i:\bigl(Y^i, \Delta_2^i\bigr)\to X$
of the MMP,
 we have $g\bigl( \Delta_1\bigr)=g^i\bigl( \Delta_1^i\bigr)$.
This is  clear for $Y^0:=Y$.
As we go from $i$ to $i+1$, the image
$g^i\bigl( \Delta_1^i\bigr)$ is unchanged if $Y^i\map Y^{i+1}$ is a flip.
Thus we need to show that
$g^{i+1}\bigl( \Delta_1^{i+1}\bigr)=g^i\bigl( \Delta_1^i\bigr)$
if $\pi_i:Y^i\to Y^{i+1}$ is a divisorial contraction with exceptional divisor
$E^i$. Let $F^i\subset E^i$ be a general fiber of $E^i\to X$.
It is clear that
$$
g^{i+1}\bigl( \Delta_1^{i+1}\bigr)\subset g^i\bigl( \Delta_1^i\bigr)
$$
and equality fails only if $E^i$ is a component of $\Delta_1^{i}$
but no other  component of $\Delta_1^{i}$ intersects $F^i$.
Since $\pi_i$ contracts a $\bigl(K_{Y^i}+\Delta_2^i\bigr)$-negative
extremal ray,
$-\Delta_1^{i} \sim_{\q,g^i}K_{Y^i}+\Delta_2^i$ shows
that $\Delta_1^{i}$ is $\pi_i$-nef. 
However, an exceptional divisor has negative intersection
with some contracted curve; a contradiction. \qed

\begin{say}[Proof of Theorem \ref{1/6.logcent.semin.thm}.1]
\label{pf.of.1/6.logcent.semin.thm} 
Set $\epsilon_i:=\mld(Z_i, X,\Delta)$. 
As in Lemma \ref{dlt.model.withD.lem}, let $g:(Y,\Delta_Y)\to (X,\Delta)$ be a
$\q$-factorial dlt model and $D_i\subset Y$ divisors such that
$a(D_i,X,\Delta)=-1+\epsilon_i$ and $ g(D_i)=Z_i$.
Set $D:=\sum_{i=1}^m D_i$; then $g(D)=Z$. 

Pick $1>c\geq 0$ such that $1-\epsilon_i\geq c$ for every $i$
and write $\Delta_Y=cD+\Delta_2$ where $\Delta_2$ is effective.
(It may have common components with $D$.)
Apply Lemma \ref{dlt.model.nefD.lem} to  get a  $\q$-factorial model
$g^m:Y^m\to X$ such that 
\begin{enumerate}
\item $\bigl(Y^m, cD^m+\Delta_2^m\bigr)$ is lc,
\item $\bigl(Y^m, \Delta_2^m\bigr)$ is dlt,
\item $K_{Y^m}+cD^m+\Delta_2^m\sim_{\q,g^m} 0$, 
\item $-D^m$ is $g^m$-nef,
\item  $\supp D^m=\supp (g^m)^{-1}(Z)$ and hence 
$g^m\bigl( D^m\bigr)=Z$.
\end{enumerate}
If  $\epsilon_i<\frac16$ for every $i$ then we can assume that $c>\frac56$.
As we noted in (\ref{5/6.first.emp}), in this case $D^m$ is seminormal and
 Lemma \ref{1/6.sn.lem} shows that $g^m_*\o_{D^m}=\o_Z$. Thus
$Z$ is   seminormal by Lemma \ref{rel.sn.lem}.\qed
\end{say}


\begin{lem} \label{1/6.sn.lem} Let $Y,X$ be normal varieties and
$g:Y\to X$ a proper morphism such that $g_*\o_Y=\o_X$.
Let $D$ be a reduced divisor on $Y$ and $\Delta''$ an
effective $\q$-divisor on $Y$. Fix some $0<c\leq 1$. Assume that
\begin{enumerate}
\item $\bigl(Y, cD+\Delta''\bigr)$ is lc,
\item $\bigl(Y, \Delta''\bigr)$ is dlt,
\item $K_Y+cD+\Delta''\sim_{\q,g}0$ and
\item $-D$ is $g$-nef (and hence $D=g^{-1}\bigl(g(D)\bigr)$.).
\end{enumerate}
Then $g_*\o_{D}=\o_{g(D)}$. 
\end{lem}

Proof. By pushing forward the exact sequence
$$
0\to \o_Y(-D)\to \o_Y\to \o_{D} \to 0
$$
we obtain
$$
\o_X=g_*\o_Y\to g_*\o_{D}\to R^1g_*\o_X(-D).
$$
Note that
$$
-D\sim_{\q,g} K_Y+\Delta''+(1-c)(-D)
$$
and the right hand side is of the form $K+\Delta+(g-\mbox{nef})$.
Let $W\subset Y$ be an lc center of $\bigl(Y,\Delta''\bigr)$.
Then $W$  is not
contained in $D$ since then  $\bigl(Y, cD+\Delta''\bigr)$
would not be lc along $W$. In particular, $D$ is  disjoint from
the general fiber of $W\to X$ by (4). Thus from Theorem \ref{amb-fuj.thm}
we conclude that none of the associated primes of $R^1g_*\o_Y(-D)$
is contained in $g(D)$. On the other hand, $g_*\o_{D}$
is supported on $g(D)$, hence 
$g_*\o_{D}\to R^1g_*\o_Y(-D)$ is the zero map.

This implies that $\o_X\to g_*\o_{D}$ is surjective. 
This map factors through $\o_{g(D)}$, hence
$g_*\o_{D}=\o_{g(D)}$. \qed

\begin{say}[A curious property of log centers] Assume that
$(X,\Delta)$ is  klt and let $Z\subset X$ be a union of  arbitrary log centers.
As in  (\ref{pf.of.1/6.logcent.semin.thm})
we construct $(Y, cD+\Delta'')$ which is klt.
Thus, as we apply  Lemma \ref{1/6.sn.lem}, the higher direct images
$R^ig_*\o_Y$ and $R^ig_*\o_Y(-D)$ are zero for $i>0$. Thus $D$ 
 is a reduced Cohen-Macaulay scheme $D$
such that 
$$
g_*\o_D=\o_{Z}\qtq{and}
R^ig_*\o_D=0 \qtq{for $i>0$.}
$$
Moreover, $D$ is a divisor on a $\q$-factorial klt pair.

This looks like a  very strong property for a reduced scheme $Z$,
but so far I have been unable to derive any useful consequences of it.
In fact, I do not know how to prove that not every reduced scheme $Z$
admits such a  morphism  $g:D\to Z$.
\end{say} 

We have used the following form of 
 \cite[3.2, 7.4]{ambro} and \cite[2.52]{fujinobook}.

\begin{thm} \label{amb-fuj.thm}
Let $g:Y\to X$ be a projective morphism and $M$ a line bundle on $Y$.
Assume that $M\sim_{\q,g} K_Y+L+\Delta$ where
$(Y,\Delta)$ is dlt and  for every log canonical center $Z\subset Y$,
the restriction of $L$ to the general fiber of $Z\to X$ is semiample.

Then every  associated prime of $R^ig_*M$ is the image of a
log canonical center of $(Y,\Delta)$.\qed
\end{thm}

\begin{say}[Proof of Theorem \ref{1/6.logcent.semin.thm}.2]
\label{pf.of.1/6.logcent.inter.thm} 

By induction on $m$, it is enough to prove Theorem \ref{1/6.logcent.semin.thm}.2
for the intersection of 2 log centers.

Let $g:(Y,\Delta_Y)\to (X,\Delta)$ be a
$\q$-factorial dlt model and $D_1, D_2\subset Y$ divisors such that
$a(D_i,X,\Delta)=-1+\mld(Z_i,X,\Delta)$ and $ g(D_i)=Z_i$.
Set $D:=D_1+D_2$. Pick any $c>0$ such that 
$\Delta_Y=cD+\Delta_2$ where $\Delta_2$ is effective
and apply Lemma \ref{dlt.model.nefD.lem}. 
Thus we get a  $\q$-factorial model
$g^m:Y^m\to X$ such that 
\begin{enumerate}
\item $K_{Y^m}+cD^m+\Delta_2^m\sim_{\q,g^m} 0$, 
\item  $\supp D^m=\supp (g^m)^{-1}(Z_1\cup Z_2)$.
\end{enumerate}
By (2), every irreducible component 
$V_j\subset Z_1\cap Z_2$ is dominated by an irreducible component of
$W_j\subset D^m_1\cap D^m_2$. By (\ref{5/6.first.emp}),
each $W_j$ is a  log center of
$\bigl(Y^m, cD^m+\Delta_2^m\bigr)$ with
$\mld\leq \mld(Z_1,X,\Delta)+\mld(Z_2,X,\Delta)$.
Thus $V_j$ is a  log center of  $(X,\Delta)$ with the same 
minimal log discrepancy.\qed
\end{say}

\begin{defn}\label{rel.semin.defn}
 Let $X$ be a reduced scheme and $U\subset X$ an open subscheme.
We say that $X$ is {\it seminormal relative to $U$}
if every finite, universal homeomorphism $\pi:X'\to X$
that is an isomorphism over $U$ is an  isomorphism.

If this holds with $U=\emptyset$, then $X$ is called {\it seminormal}.
For more details, see \cite[Sec.I.7.2]{rc-book}.

If $X$ satisfies Serre's condition $S_2$ then 
seminormality depends only on the codimension 1 points of $X$.
That is,  $X$ is  seminormal relative to $U$
iff there is a closed subset $Z\subset X$ of codimension $\geq 2$
such that $X\setminus Z$ is  seminormal relative to $U$.
\end{defn}

With this definition, we can state the theorem behind 
Corollary  \ref{1/2.coeff.seminorm.thm} as follows.

\begin{thm} \label{1/2.logcent.semin.thm}
Let $(X,S+\Delta)$ be an lc pair where $S$ is $\q$-Cartier.
Let $Z_i\subset X$ be  log centers of $(X,\Delta)$
 for  $i=1,\dots, m$.

 If  $\mld(Z_i,X,\Delta)<\frac12$ for every $i$ then
 $S\cup Z_1\cup\cdots\cup Z_m$ is
 seminormal relative to $X\setminus S$.
\end{thm}

Proof. By passing to a cyclic cover 
and  using  Lemma \ref{rel.sn.desc.lem}   we may assume that
$S$ is Cartier.

Next we closely follow 
(\ref{pf.of.1/6.logcent.semin.thm}).  
Let $g:(Y,S_Y+\Delta_Y)\to (X,S+\Delta)$ be a
$\q$-factorial dlt model and $D_i\subset Y$ divisors such that
$a(D_i,X,\Delta)=-1+\mld(Z_i,X,\Delta)$ and $ g(D_i)=Z_i$.
Pick $c>\frac12$ such that $1-\mld(Z_i,X,\Delta)\geq c$ for every $i$.
Set $D:=S_Y+\sum_i D_i$
and write $\Delta_Y=cD+\Delta_2$ where $\Delta_2$ is effective.

Apply Lemma \ref{dlt.model.nefD.lem}  to get a  $\q$-factorial model
$g^m:Y^m\to X$ such that 
\begin{enumerate}
\item $\bigl(Y^m, cD^m+\Delta_2^m\bigr)$ is lc,
\item $\bigl(Y^m,  \Delta_2^m\bigr)$ is dlt,
\item $K_{Y^m}+ cD^m+\Delta_2^m\sim_{\q,g^m} 0$, 
\item $-D^m$ is $g^m$-nef and 
\item   $g^m\bigl( D^m\bigr)=S\cup Z_1\cup\cdots\cup Z_m$.
\end{enumerate}

Using Lemmas \ref{1/6.sn.lem} and \ref{rel.sn.lem}
we see that it is enough to prove that
$D^m$ is seminormal relative to $Y^m\setminus S_Y^m$.

Since $Y^m$ is dlt, it is Cohen--Macaulay hence
$D^m$ is $S_2$. As  we noted in 
Definition \ref{rel.semin.defn}, it is sufficient to check seminormality at
codimension 2 points of $Y^m$. 
As in (\ref{5/6.first.emp}), this reduces to the smooth surface case.
We see that if $F$ is a smooth surface,
 $(F, S+cD)$ is lc  and $c>\frac12$ then, 
at every point of
$S\cap D$, $D$ is smooth and intersects $S$ transversally.
Thus $S+D$ is seminormal at all points of $S\cap D$.
\qed
\medskip

Again note that the bound $\frac12$ is sharp;
  $\bigl(\a^2, (x=0)+\frac12(x+y=0)+\frac12(x-y=0)\bigr)$
is lc  but its boundary is not  seminormal at the origin. 

\begin{say}[Proof of Corollary \ref{1/2.coeff.seminorm.thm}] 
None of the irreducible components of $\Delta$ is contained in a
fiber of $f$, hence $f:B_J\to C$ is flat. The main point is
to show that its fibers are reduced. 

If $b_j>\frac12$ then the corresponding divisor $B_j$
is a log center and
$\mld(B_j,X,\Delta)=1-b_j<\frac12$. 
Thus, by Theorem \ref{1/2.logcent.semin.thm},
$X_c+B_J$ is seminormal relative to $X\setminus X_c$ for every $c\in C$.
By Lemma \ref{sn.coms.inters.red} this implies that
$X_c\cap B_J$ is reduced. \qed
\end{say}

We have used three easy properties of seminormal schemes.

\begin{lem} \label{rel.sn.desc.lem}
Let $g:Y\to X$ be a finite morphism of normal schemes.
Let $Z\subset X$ be a closed, reduced subscheme and 
 $U\subset X$ an open subscheme.
If  $\red g^{-1}(Z)$ is  seminormal relative to $g^{-1}U$ then
 $Z$ is  seminormal relative to $U$.
\end{lem}

Proof. We may assume that $X,Y$ are irreducible and affine.
Let $\pi:Z'\to Z$ be a  finite, universal homeomorphism 
that is an isomorphism over $Z\cap U$.
Pick $\phi\in \o_{Z'}$. Since  $\red g^{-1}(Z)$ is  seminormal relative to 
$g^{-1}U$, the pull back
$\phi\circ g$ is a regular function on $\red g^{-1}(Z)$.
 We can lift it  to a regular function
$\Phi_X$ on $X$. Since $Y$ is normal,
$$
\Phi_Y:=\tfrac1{\deg X/Y}\tr_{X/Y} \Phi_X
$$
is  regular on $Y$ and  $\Phi_Y|_{Z}=\phi$. 
Thus  $Z$ is  seminormal relative to $U$.\qed 

\begin{lem} \label{rel.sn.lem}
Let $g:Y\to X$ be a proper morphism of reduced schemes
such that $g_*\o_Y=\o_X$. Let $U\subset X$ be an open subscheme.
If  $Y$ is  seminormal relative to $g^{-1}U$ then
 $X$ is  seminormal relative to $U$.
\end{lem}

Proof. Let $\pi:X'\to X$ be a  finite, universal homeomorphism 
that is an isomorphism over $U$.
Set $Y':=\red\bigl(Y\times_XX')$ with projection
$\pi_Y:Y'\to Y$. Then $\pi_Y$ is a  finite, universal homeomorphism 
that is an isomorphism over $g^{-1}U$. Thus $\pi_Y$ is an isomorphism,
so we can factor $g$ as  $Y\to X'\to X$.
This implies that $\pi_*\o_{X'}\subset g_*\o_Y=\o_X$,
hence $\pi$  is an  isomorphism. \qed

\begin{lem} \label{sn.coms.inters.red} Let $X$ be semi normal relative to $U$.
Let $X_1, X_2\subset X$ be closed, reduced subschemes such that
$X=X_1\cup X_2$.
 Then $\o_{X_1\cap X_2}$ has no nilpotent elements whose support is in
$X\setminus U$.
\end{lem}

Proof. Let $I\subset \o_{X_1\cap X_2}$ be the ideal sheaf of
 nilpotent elements whose support is in
$X\setminus U$ and $r(X_1\cap X_2)\subset X_1\cap X_2$
the corresponding subscheme.

Let $r_i:\o_{X_i}\to \o_{X_1\cap X_2}$ 
and $\bar r_i:\o_{X_i}\to \o_{r(X_1\cap X_2)}$ denote the
restriction maps. Then $\o_X$ sits in an exact sequence
$$
0\to \o_X\to \o_{X_1}+\o_{X_2}\stackrel{(r_1, -r_2)}{\longrightarrow}
 \o_{X_1\cap X_2}\to 0.
$$
The similar sequence
$$
0\to A\to \o_{X_1}+\o_{X_2}\stackrel{(\bar r_1, -\bar r_2)}{\longrightarrow}
 \o_{r(X_1\cap X_2)}\to 0
$$
defines a coherent sheaf of $\o_X$-algebras $A$
and $\spec_XA\to X$ is a 
 finite, universal homeomorphism $\pi:X'\to X$
that is an isomorphism over $U$. 
Since $X$ is  semi normal relative to $U$, $A=\o_X$ hence
$X_1\cap X_2=r(X_1\cap X_2)$. \qed

\begin{ack} I thank V.~Alexeev and O.~Fujino 
 for useful  comments and F.~Ambro for many corrections and
remarks.
Partial financial support   was provided by  the NSF under grant number 
DMS-0758275.
\end{ack}

\newcommand{\etalchar}[1]{$^{#1}$}
\providecommand{\bysame}{\leavevmode\hbox to3em{\hrulefill}\thinspace}
\providecommand{\MR}{\relax\ifhmode\unskip\space\fi MR }
\providecommand{\MRhref}[2]{%
  \href{http://www.ams.org/mathscinet-getitem?mr=#1}{#2}
}
\providecommand{\href}[2]{#2}

\bigskip

\noindent Princeton University, Princeton NJ 08544-1000, USA

\begin{verbatim}kollar@math.princeton.edu\end{verbatim}

\end{document}